\documentclass[11pt]{article}
\usepackage{amssymb,amsmath,latexsym,amsthm}
\usepackage{url}
\usepackage{enumerate}
\usepackage[bookmarks=false,colorlinks=true,linkcolor={blue},pdfstartview={XYZ null null 1.22}]{hyperref}
\usepackage{mathrsfs,eufrak}
\usepackage{units}

\usepackage[dvips]{graphicx}

\usepackage{pstricks}

\newtheorem{prop}{Proposition}[section]

\newtheorem{cor}[prop]{Corollary}

\newtheorem{lemma}[prop]{Lemma}

\newtheorem{thm}[prop]{Theorem}
\newtheorem{defi}[prop]{Definition}

\renewcommand{\geq}{\geqslant}
\def\leq{\leqslant}

\newcommand{\Z}{\mathbb{Z}}
\newcommand{\R}{\mathbb{R}}

\newcommand{\ff}{\varphi}

\newcommand{\poly}{Q}

\def\e{\varepsilon}

\def\1{{\mathbf{1}}}

\def\1{{\mathbf{1}}}
\def\0.5{{\frac{1}{2}}}

\title{Central and non-central limit theorems\\ in a free probability setting}

\author{
  Ivan Nourdin\thanks{Ivan Nourdin was partially supported by the ANR Grants ANR-09-BLAN-0114 and ANR-10-BLAN-0121
at Universit\'e de Lorraine.}\\
\small
  Universit\'e de Lorraine,
  Institut de Math\'ematiques \'Elie Cartan\\
\small
  Facult\'e des Sciences et Techniques, Campus Aiguillettes,\\
\small
  B.P. 70239,
  54506 Vandoeuvre-l\`es-Nancy Cedex, France\\
\and
\small
  and\\
\small
  Fondation des Sciences Math\'ematiques de Paris\\
\small
  IHP, 11 rue Pierre et Marie Curie,
  75231 Paris Cedex 05, France\\
\small
 {\tt inourdin@gmail.com}\\
\\
\and
  Murad S. Taqqu\thanks{Murad S. Taqqu was partially supported by the NSF Grant DMS-1007616 at Boston University.}\\
\small
  Boston University, Departement of Mathematics\\
\small
  111 Cummington Road, Boston (MA), USA\\
\small
  {\tt murad@math.bu.edu}  \\
  }

\begin{document}

\maketitle

\begin{abstract}

Long-range dependence in time series may yield non-central limit theorems.
We show that there are analogous time series in free probability with limits
represented by multiple Wigner integrals, where Hermite processes are replaced
by non-commutative Tchebycheff processes. This includes the non-commutative fractional Brownian motion
and the non-commutative Rosenblatt process.
\bigskip


\bigskip


\end{abstract}

\eject

\section{Introduction and main result}\label{s:intro}

Normalized sums of i.i.d.\!\! random variables satisfy the usual central limit theorem.
But it is now well-known that this is not necessarily the case if the i.i.d.\!\! random
variables are replaced by a stationary sequence with long-range dependence, that is,
with a correlation which decays slowly as the lag tends to infinity.
We want to investigate whether similar {\em non-central} results hold in the free
probability setting.

We are motivated by the fact that there is often a close
correspondence between classical probability and free probability. For example, the
Gaussian law has the semicircular law as an analogue, hence the notion of a
stationary {\em semicircular sequence} with a given correlation function. Multiple
Wiener integrals which span the so-called Wiener chaos have {\em multiple Wigner
integrals} as an analogue, with Hermite polynomials being replaced by {\em Tchebycheff
polynomials}. We will see that the notion of Hermite rank is to be replaced by the
{\em Tchebycheff rank} and that long-range dependence in the free probability
setting also yields non-standard limits which are somewhat analogous to those
in classical probability. In classical probability, the limits can be represented
by multiple Wiener integrals. In free probability, they are represented by multiple
Wigner integrals with similar kernels. Finally, the Hermite processes that
appear in the limit in the usual probability setting are replaced, in the
free probability setting, by {\em non-commutative Tchebycheff processes}.

Before stating our main results,
let us describe their analogue in the classical probability framework.
Let $Y=\{Y_k:\,k\in\Z\}$ be a stationary {\em Gaussian} sequence on a probability space $(\Omega,\mathcal{F},P)$, with $E[Y_k]=0$ and $E[Y_k^2]=1$,
and let $\rho(k-l)=E[Y_kY_l]$ be its correlation kernel.
(Observe that $\rho$ is symmetric, that is, $\rho(n)=\rho(-n)$ for all $n\geq 1$.)
Let $H_0(y)=1$, $H_1(y)=y$, $H_2(y)=y^2-1$, $H_3(y)=y^3-3y$, $\ldots$, denote the
sequence of Hermite polynomials
(determined by the recursion $yH_k=H_{k+1}+kH_{k-1}$), and consider a real-valued
polynomial $Q$ of the form
\begin{equation}\label{decompo-cl}
Q(y)=\sum_{s\geq q} a_s H_s(y),
\end{equation}
with $q\geq 1$ and $a_q\neq 0$, and where only a finite number of coefficients $a_s$ are
non zero.
The integer $q$ is called the {\em Hermite rank} of $Q$.
Finally, set
\begin{equation}\label{vnpt-cl}
W_n(Q,t)=\sum_{k=1}^{[nt]} Q(Y_k)=\sum_{s\geq q}a_s\sum_{k=1}^{[nt]}H_s(Y_k),
\quad t\geq 0.
\end{equation}

The following theorem is a summary of the main findings in Breuer and Major \cite{breuermajor},
Dobrushin and Major \cite{dobrushin-major} and Taqqu \cite{Taqqu1,Taqqu2}.
Here and throughout the sequel, the notation `f.d.d.' stands for the convergence in finite-dimensional distribution.
\begin{thm}\label{main-cl}
Let $Q$ be the polynomial defined by (\ref{decompo-cl}) and let $W_n(Q,\cdot)$ be defined by (\ref{vnpt-cl}).
\begin{enumerate}
\item If $\sum_{k\in\Z}|\rho(k)|^q$ is finite,
then, as $n\to\infty$,
\begin{equation}\label{tcl1-cl}
\frac{W_n(Q,\cdot)}{\sqrt{n}}\,\overset{\rm f.d.d.}{\to}\,\sqrt{\sum_{s\geq q}s!
a_s^2\sum_{k\in\Z}\rho(k)^s}\times B,
\end{equation}
with $B$ a classical Brownian motion.
\item Let $L:(0,\infty)\to (0,\infty)$ be a function which is slowly varying at infinity and
bounded away from 0 and infinity on every compact subset of $[0,\infty)$.
If $\rho$ has the form
\begin{equation}\label{rho}
\rho(k)= k^{-D}L(k),\quad k\geq 1,
\end{equation}
with $0<D<\frac1q$,
then, as $n\to\infty$,
\begin{equation}\label{tcl2-cl}
\frac{W_n(Q,\cdot)}{n^{1-qD/2}L(n)^{q/2}}\,\overset{\rm f.d.d.}{\to}\,
\frac{a_q}{\sqrt{(1-\frac{qD}2)(1-qD)}}\times R_{1-qD/2,q},
\end{equation}
where $R_{1-qD/2,q}$ is the $q$th Hermite process of parameter $1-qD/2$.
\end{enumerate}
\end{thm}

For background on the notions of Hermite
processes, see
Dobrushin and Major \cite{dobrushin-major},
Embrechts and Maejima \cite{embrechts-maejima},
Peccati and Taqqu \cite{peccati-taqqu:2011-book}, and
Taqqu \cite{Taqqu1,Taqqu2,taqqu:2003-livre}.
A slowly varying function at infinity $L(x)$ is such that $\lim_{x\to\infty} L(cx)/L(x)=1$ for all $c>0$.
Constants and logarithm are slowly varying.
A useful property is the {\em Potter's bound} (see \cite[Theorem 1.5.6, $(ii)$]{BGT}): for every $\delta>0$,
there is $C=C(\delta)>1$ such that, for all $x,y>0$,
\begin{equation}\label{potter}
\frac{L(x)}{L(y)}\leq C\,\max\left\{ \left(\frac{x}{y}\right)^{\delta},
\left(\frac{x}{y}\right)^{-\delta}
\right\}.
\end{equation}

When $\rho$ is given by (\ref{rho}), then
\[
\sum_{k\in\Z}|\rho(k)|^q
\left\{
\begin{array}{ll}
<\infty&\quad\mbox{if $D>1/q$}\\
=\infty&\quad\mbox{if $0<D<1/q$}\\
\end{array}
\right..
\]
When $D>1/q$, one says that the process $X$ has {\em short-range dependence}.
When $0<D<1/q$, one says that it has {\em long-range dependence}.
In the critical case $D=1/q$, the series may be finite or infinite, depending on the precise
value of $L$; this is the reason why we do not investigate this case further in Theorem \ref{main-cl} nor in
the forthcoming Theorem \ref{main}.

We now state our main result, which may be regarded as a non-commutative counterpart of Theorem \ref{main-cl}.
Some of the terms and concepts of free probability used here
are defined in Section \ref{prelimi} which introduces free probability
in a nutshell.

Let $X=\{X_k:\,k\in\Z\}$ be a stationary semicircular sequence on a non-commutative probability space $(\mathscr{A},\ff)$,
assume that $\ff(X_k)=0$ and $\ff(X_k^2)=1$
and let $\rho(k-l)=\ff(X_kX_l)$ be its correlation kernel.
(Observe that $\rho$ is symmetric, that is, $\rho(n)=\rho(-n)$ for all $n\geq 1$.)
Let $U_0(x)=1$, $U_1(x)=x$, $U_2(x)=x^2-1$, $U_3(x)=x^3-2x$, $\ldots$, denote the
sequence of Tchebycheff polynomials of second kind
(determined by the recursion $xU_k=U_{k+1}+U_{k-1}$). 
For a presentation of polynomials in the classical and free probability setting, see Anshelevich \cite{ans}.
Consider a real-valued
polynomial $Q$ of the form
\begin{equation}\label{decompo}
Q(x)=\sum_{s\geq q} a_s U_s(x),
\end{equation}
with $q\geq 1$ and $a_q\neq 0$, and where only a finite number of coefficients $a_s$ are
non zero.
The integer $q$ is called the {\em Tchebycheff rank} of $Q$.
Finally, set
\begin{equation}\label{vnpt}
V_n(Q,t)=\sum_{k=1}^{[nt]} Q(X_k)=\sum_{s\geq q}a_s\sum_{k=1}^{[nt]}U_s(X_k),
\quad t\geq 0.
\end{equation}

Our main result goes as follows:
\begin{thm}\label{main}
Let $Q$ be the polynomial defined by (\ref{decompo}) and let $V_n(Q,\cdot)$ be defined by (\ref{vnpt}).
\begin{enumerate}
\item If $\sum_{k\in\Z}|\rho(k)|^q$ is finite,
then, as $n\to\infty$,
\begin{equation}\label{tcl1}
\frac{V_n(Q,\cdot)}{\sqrt{n}}\,\overset{\rm f.d.d.}{\to}\,\sqrt{\sum_{s\geq q}
a_s^2\sum_{k\in\Z}\rho(k)^s}\times S,
\end{equation}
with $S$ a free Brownian motion, defined in Section \ref{section fBM}.
\item Let $L:(0,\infty)\to (0,\infty)$ be a function which is slowly varying at infinity and
bounded away from 0 and infinity on every compact subset of $[0,\infty)$.
If $\rho$ has the form (\ref{rho})
with $0<D<\frac1q$,
then, as $n\to\infty$,
\begin{equation}\label{tcl2}
\frac{V_n(Q,\cdot)}{n^{1-qD/2}L(n)^{q/2}}\,\overset{\rm f.d.d.}{\to}\,
\frac{a_q}{\sqrt{(1-\frac{qD}2)(1-qD)}}\times R_{1-qD/2,q},
\end{equation}
where $R_{1-qD/2,q}$ is the $q$th non-commutative Tchebycheff process of parameter $1-qD/2$,
defined by a multiple Wigner integral of order $q$. It is
given in Definition \ref{d:tcheb}.
\end{enumerate}
\end{thm}

Let us compare Theorem \ref{main-cl} and Theorem \ref{main}.
The fact that (\ref{tcl1}) relies on the free Brownian motion $S$ implies for example that the marginal distribution of $S$ is not Gaussian but is
the semicircular law defined in Section \ref{semi-sec}, which has, in particular, a compact support.
The notion of free independence is very different from the classical
notion of independence as noted in Section \ref{free-sec}.
The resulting multiple integrals are then of a very different nature.
Also, since the Hermite and Tchebycheff polynomials
are different, the decomposition of the polynomial $Q$ in Tchebycheff polynomials is different from
its decomposition in Hermite
polynomials and, consequently, its Tchebycheff rank can be different from its Hermite rank.
This implies that even the order of the multiple integral in the limit may not be the same.

In the non-commutative probability context,
the exact expression of
the underlying law is of relatively minor importance. Moments are essential. In fact, convergence in law
can be defined in our setting
through convergence of moments.
Also, in the classical probability setting, one can consider a function $Q(x)$ in (\ref{vnpt-cl}) which is not
necessarily a polynomial. One could possibly do so as well in the non-commutative probability setting.

The methods of proof of Theorems \ref{main-cl} and \ref{main} are different. The proof of the central limit theorem
in Theorem \ref{main} does not need to rely on cumulants and diagrams as in the classical case (see \cite{breuermajor}). It uses
a {\em transfer principle} established recently in \cite{NPS} and stated in Proposition \ref{transfer}. The proof of the non-central limit theorem
in Theorem \ref{main} is much simpler than in the classical proof of Dobrushin and Major \cite{dobrushin-major}
and Taqqu \cite{Taqqu1,Taqqu2} because it is sufficient here to establish the convergence of joint moments.
Such an approach breaks down in the classical case when $q\geq 3$ because, there, the joint moments do not
characterize the target distribution anymore.

We proceed as follows. In Section \ref{prelimi}, we present free probability in a nutshell and define multiple
Wigner integrals. In Section \ref{s:ncTp}, we introduce the non-commutative fractional Brownian motion, the
non-commutative Rosenblatt process and more generally the non-commutative Tchebycheff processes, and study
some of their  basic properties. We compute their joint moments in Section \ref{s:joint}.
Theorem \ref{main} is proved in Section  \ref{s:proof}.

\section{Free probability in a nutshell}\label{prelimi}

\subsection{Random matrices}

Let $(X_{ij})_{i,j\geq 1}$ be a sequence of independent standard Brownian motions, all defined on the same probability space
$(\Omega,\mathcal{F},P)$. Consider the random matrices
\[
A_n(t)=\left(\frac{X_{ij}(t)}{\sqrt{n}}\right)_{1\leq i,j\leq n},
\]
and
\begin{equation}\label{mn}
M_n(t)=\frac{A_n(t)+A_n(t)^*}{\sqrt{2}},
\end{equation}
where $A_n(t)^*$ denotes the transpose or adjoint of the matrix $A_n(t)$. Thus, $M_n(t)$ is self-adjoint.
For each $t$, $A_n(t)$ and $M_n(t)$ both belong to $\mathscr{A}_n$, the set of random matrices with entries in $L^{\infty -}(\Omega)=\cap_{p>1} L^p(\Omega)$,
that is, with all moments.
On $\mathscr{A}_n$, consider the linear form $\tau_n:\mathscr{A}_n\to\R$ defined by
\begin{equation}\label{e:trace}
\tau_n(M)=E\left[\frac{1}n{\rm Tr}(M)\right],
\end{equation}
where $E$ denotes the mathematical
expectation associated to $P$, whereas ${\rm Tr}(\cdot)$ stands for the usual trace operator.
The space $(\mathscr{A}_n,\tau_n)$ is the prototype of a {\em non-commutative} probability space.

Let $t_k>\ldots>t_1>t_0=0$. A celebrated theorem by Voiculescu \cite{Voiculescu2} asserts that
the increments $M_n(t_1),M_n(t_2)-M_n(t_1),\ldots,M_n(t_k)-M_n(t_{k-1})$
are asymptotically free, meaning that
\[
\tau_n\left(
Q_1(M_n(t_{i_1})-M_n(t_{i_1-1}))\ldots Q_m(M_n(t_{i_m})-M_n(t_{i_m-1}))
\right)\to 0\quad\mbox{as $n\to\infty$},
\]
for all $m\geq 2$, all $i_1,\ldots,i_m\in\{1,\ldots,k\}$ with
$i_1\neq i_2$, $i_2\neq i_3$, $\ldots$, $i_{m-1}\neq i_m$, and all
real-valued polynomials $Q_1,\ldots,Q_m$ such that $\tau_n(Q_i(M_n(t_i)-M_n(t_{i-1})))\to 0$ as $n\to\infty$ for each $i=1,\ldots,m$.

Let $t>0$. The celebrated Wigner theorem \cite{wigner}
can be formulated as follows: as $n\to\infty$,
$M_n(t)$ converges in law to the semicircular law of variance $t$, that is, for any real-valued polynomial $Q$,
\[
\tau_n(Q(M_n(t)))\to\frac{1}{2\pi t}\int_{-2\sqrt{t}}^{2\sqrt{t}} Q(x)\sqrt{4t-x^2}dx.
\]

In the same way as calculus provides a nice setting for studying limits of sums and as classical Brownian motion
provides a nice setting for studying limits of random walks, free probability provides a convenient framework for
investigating limits of random matrices. In the setting of free probability, free independence (see Section \ref{free-sec}) replaces independence,
increments have a semicircular marginal law (see Section \ref{semi-sec}) and thus, free Brownian motion
(see Section \ref{section fBM}) which has these properties,
can be used to study random matrices $M_n(t)$ for large $n$. Conversely, one may visualize free Brownian
motion as a large random matrix $M_n(t)$.

\subsection{Non-commutative probability space} \label{prelimi-free}
In this paper, we use the phrase ``{\em non-commutative probability space}''
to indicate a {\em von Neumann algebra} $\mathscr{A}$ (that is,
an algebra of operators on a complex
separable Hilbert space,
closed under adjoint and convergence in the weak operator topology) equipped with a {\em trace} $\ff$,
that is, a {\em unital linear functional} (meaning preserving the identity) which is weakly
continuous, positive (meaning $\ff(X)\ge 0$
whenever $X$ is a non-negative element of $\mathscr{A}$; i.e.\ whenever $X=YY^\ast$
for some $Y\in\mathscr{A}$), {\em faithful} (meaning that if
$\ff(YY^\ast)=0$ then $Y=0$), and {\em tracial} (meaning that $\ff(XY)=\ff(YX)$ for all
$X,Y\in\mathscr{A}$, even though in general $XY\ne YX$).
We will not need to use the full force of this definition, only some of its
consequences. See \cite{nicaspeicher} for a systematic presentation.

\subsection{Random variables}\label{s:rv}
In a non-commutative probability space, we refer to the self-adjoint elements of the
algebra as
{\em random variables}.  Any
random variable $X$ has
a {\em law}: this is the unique
probability measure $\mu$ on $\R$
with the same moments as $X$; in other words, $\mu$ is such that
\begin{equation}\label{mu}
\int_{\R} \poly(x) d\mu(x) = \ff(\poly(X)),
\end{equation}
for any polynomial $\poly$.
(The existence and uniqueness of $\mu$ follow from the positivity of $\ff$,  see \cite[Proposition 3.13]{nicaspeicher}.)
Thus $\varphi$ acts as an expectation. The $\tau_n$ in (\ref{e:trace}), for example,
play the role of $\varphi$.
Also, while there is a classical random variable $\widehat{X}\in\R$ with law
$\mu$, there is, a priori, no direct relationship between $X$ and $\widehat{X}$.

\subsection{Convergence in law}
We say that a sequence $(X_{1,n},\ldots,X_{k,n})$, $n\geq 1$, of random vectors {\em converges in law}
to a random vector $(X_{1,\infty},\ldots,X_{k,\infty})$, and we write \[
(X_{1,n},\ldots,X_{k,n})\overset{\rm law}{\to} (X_{1,\infty},\ldots,X_{k,\infty}),
\]
to indicate the convergence in the sense of (joint) moments, that is,
\begin{equation}\label{cvlaw}
\lim_{n\to\infty}\ff\left(\poly(X_{1,n},\ldots,X_{k,n})\right) = \ff\left(\poly(X_{1,\infty},\ldots,X_{k,\infty})\right),
\end{equation}
for any  polynomial $\poly$ in $k$ non-commuting variables.
In the case of vectors, there may be no corresponding probability law $\mu_{(X_1,\ldots,X_n)}$ as in (\ref{mu}),
see \cite[Lecture 4]{nicaspeicher}.

We say that a sequence $(F_n)$ of {\em non-commutative stochastic processes} (that is, each $F_n$ is a
one-parameter family of
self-adjoint operators $F_n(t)$ in the non-commutative probability space $(\mathscr{A},\ff)$) {\em converges in the sense of finite-dimensional
distributions} to a non-commutative stochastic process $F_\infty$,
and we write \[
F_{n}\overset{\rm f.d.d.}{\to} F_\infty,
\]
to indicate that, for any $k\geq 1$ and any $t_1,\ldots,t_k\geq 0$,
\[
(F_{n}(t_1),\ldots,F_{n}(t_k))\overset{\rm law}{\to} (F_{\infty}(t_1),\ldots,F_{\infty}(t_k)).
\]

\subsection{Free independence}\label{free-sec}

In the non-commutative probability setting, the notion of {\em independence} (introduced by
Voiculescu in \cite{Voiculescu})
goes as follows. Let $\mathscr{A}_1,\ldots,\mathscr{A}_p$ be unital subalgebras of $\mathscr{A}$.  Let $X_1,\ldots, X_m$ be elements
chosen from among the $\mathscr{A}_i$'s such that, for $1\le j<m$,
two consecutive elements $X_j$ and $X_{j+1}$ do not come from the same $\mathscr{A}_i$, and
such that $\ff(X_j)=0$ for each $j$.  The subalgebras $\mathscr{A}_1,\ldots,\mathscr{A}_p$ are said to be {\em free} or {\em freely
independent} if, in this circumstance,
\begin{equation}\label{free-def}
\ff(X_1X_2\cdots X_m) = 0.
\end{equation}
Random variables are called freely independent if the unital algebras
they generate are freely independent.  Freeness is in general much more complicated than classical independence.
Nevertheless, if $X,Y$ are freely independent, then their joint moments
are determined by the moments of $X$ and $Y$ separately
as in the classical case.
For example, if
$X,Y$ are free and $m,n\geq 1$, then by (\ref{free-def}),
\[
\ff\big((X^m-\ff(X^m)1)(Y^n-\ff(Y^n)1)\big)=0.
\]
By expanding (and using the linear property of $\ff$), we get
\begin{equation}\label{covcov}
\ff(X^mY^n)=\ff(X^m)\ff(Y^n),
\end{equation}
which is what we would expect under classical independence. But, by setting
$X_1=X_3=X-\ff(X)1$ and $X_2=X_4=Y-\ff(Y)$
in (\ref{free-def}), we note that two consecutive $X_j$ do not belong to the same subalgebra and hence, by (\ref{free-def}), we also have
\[
\ff\big((X-\varphi(X)1)(Y-\varphi(Y)1)(X-\varphi(X)1)(Y-\varphi(Y)1)\big)=0.
\]
By expanding, using (\ref{covcov}) and the tracial property of $\ff$ (for instance $\ff(XYX)=\ff(X^2Y)$)
we get
\begin{eqnarray*}
\ff(XYXY)
&=&\ff(Y)^2\ff(X^2)+\ff(X)^2\ff(Y^2)-\ff(X)^2\ff(Y)^2,
\end{eqnarray*}
which is different from $\ff(X^2)\ff(Y^2)$, which is what one would have obtained if $X$ and $Y$ were
classical independent random variables.
Let us note, furthermore, that the relation between moments and cumulants\footnote{
Cumulants have the following property which linearizes independence:
\begin{equation}\label{cum}
\kappa_n(X+Y,\ldots,X+Y)=\kappa_n(X,\ldots,X)+\kappa_n(Y,\ldots,Y),\quad n\geq 1.
\end{equation}
Relation (\ref{cum}) holds in classical probability if $X$ and $Y$ are independent
random variables and it holds in free probability if $X$ and $Y$ are freely independent
(see \cite[Proposition 12.3]{nicaspeicher}). Since the classical notion of independence is different
from the notion of free independence, the cumulants $\kappa_n$ in classical probability
are different from those in free probability.
} is different from the classical case
(see \cite[identity (11.8)]{nicaspeicher}).

\subsection{Semicircular distribution}\label{semi-sec}

The {\em semicircular distribution}
$\mathcal{S}(m,\sigma^2)$
with mean $m\in\R$ and variance $\sigma^2>0$
 is the probability distribution
\begin{equation} \label{eq semicircle}
\mathcal{S}(m,\sigma^2)(dx) = \frac{1}{2\pi \sigma^2} \sqrt{4\sigma^2-(x-m)^2}\,{\bf 1}_{\{|x-m|\le 2\sigma\}}\,dx.
\end{equation}
If $m=0$, this distribution is symmetric around $0$,
and therefore its odd moments are all $0$. A simple calculation shows that the even centered moments are given by
(scaled) {\em Catalan numbers}: for non-negative integers $k$,
\[ \int_{m-2\sigma}^{m+2\sigma} (x-m)^{2k} \mathcal{S}(m,\sigma^2)(dx) = C_k \sigma^{2k}, \]
where
\[
C_k = \frac{1}{k+1}\binom{2k}{k}
\]
(see, e.g., \cite[Lecture 2]{nicaspeicher}).
In particular, the variance is $\sigma^2$ while the centered fourth moment is $2\sigma^4$.
The semicircular distribution plays here the role of the Gaussian distribution.
It has the following similar properties:
\begin{enumerate}
\item If $S\sim\mathcal{S}(m,\sigma^2)$ and $a,b\in\R$, then $aS+b\sim \mathcal{S}(am+b,a^2\sigma^2)$.
\item If $S_1\sim\mathcal{S}(m_1,\sigma_1^2)$ and $S_2\sim\mathcal{S}(m_2,\sigma_2^2)$ are freely independent,
then $S_1+S_2\sim\mathcal{S}(m_1+m_2,\sigma_1^2+\sigma_2^2)$.
\end{enumerate}
The second property can be readily verified using the $R$-transform. The $R$-transform of a random variable
$X$ is the generating function of its free cumulants. It is such that, if $X$ and $Y$ are freely independent, then
$R_{X+Y}(z)=R_{X}(z)+R_{Y}(z)$. The $R$-transform of the semicircular law is
$R_{\mathcal{S}(m,\sigma^2)}(z)=m+\sigma^2z$, $z\in\mathbb{C}$ (see \cite[Formula (11.13)]{nicaspeicher}).

\subsection{Free Brownian Motion} \label{section fBM}

A {\em one-sided free Brownian motion} $S=\{S(t)\}_{t\geq 0}$ is a non-commutative stochastic process
with the following defining characteristics:
\begin{itemize}
\item[(1)] $S(0) = 0$.
\item[(2)] For $t_2>t_1\geq 0$, the law of $S(t_2)-S(t_1)$ is the semicircular distribution of mean 0 and variance $t_2-t_1$.
\item[(3)] For all $n$ and $t_n>\cdots>t_2>t_1>0$, the increments $S(t_1)$, $S(t_2)-S(t_1)$, \ldots, $S(t_n)-S(t_{n-1})$ are
freely independent.
\end{itemize}
A {\em two-sided free Brownian motion} $S=\{S(t)\}_{t\in\R}$ is defined to be
\[
S(t)=\left\{
\begin{array}{ll}
S_1(t)&\quad\mbox{if $t\geq 0$}\\
S_2(-t)&\quad\mbox{if $t<0$}
\end{array}
\right.,
\]
where $S_1$ and $S_2$ are two freely independent one-sided free Brownian motions.

\subsection{Wigner integral}\label{wigner1st}

From now on, we suppose that $L^2(\R^p)$
stands for the set of all {\em real-valued} square-integrable functions on $\R^p$. When $p=1$, we only write $L^2(\R)$ to
simplify the notation.

Let $S=\{S(t)\}_{t\in\R}$ be a two-sided free Brownian motion.
Let us quickly sketch out the construction of the {\em Wigner integral} of $f$ with respect to $S$.
For an indicator function $f={\bf 1}_{[u,v]}$, the Wigner integral of $f$ is defined
by
\[
\int_\R {\bf 1}_{[u,v]}(x)dS(x)=S(v)-S(u).
\]
We then extend this definition by linearity to simple functions of the form
$
f=\sum_{i=1}^k \alpha_i {\bf 1}_{[u_i,v_i]},
$
where $[u_i,v_i]$ are disjoint intervals of $\R$.
Simple computations show that
\begin{eqnarray}
\ff\left(\int_\R f(x)dS(x)\right)&=&0\label{wigner1}\\
\ff\left(\int_\R f(x)dS(x)\times\int_\R g(x)dS(x)\right)&=&\langle f,g\rangle_{L^2(\R)}.\label{wigner2}
\end{eqnarray}
By approximation, the definition of $\int_\R f(x)dS(x)$ is extended to all $f\in L^2(\R)$,
and (\ref{wigner1})-(\ref{wigner2}) continue to hold in this more general setting.
As anticipated, the Wigner process $\left\{\int_\R f(x)dS(x):\,f\in L^2(\R)\right\}$ forms a
centered semicircular family in the sense of the forthcoming Section \ref{sec:semfree}.

\subsection{Semicircular sequence and semicircular process}\label{sec:semfree}

Let $k\geq 2$. A random vector $(X_1,\ldots,X_k)$ is said to have a {\em $k$-dimensional
semicircular distribution} if, for every $\lambda_1,\ldots,\lambda_k\in\R$, the random variable
$\lambda_1 X_1 +\ldots + \lambda_k X_k$ has a semicircular distribution.
In this case, one says that the random variables $X_1,\ldots,X_k$ are {\em jointly semicircular} or,
alternatively, that $(X_1,\ldots,X_k)$ is a {\em semicircular vector}.
As an example, one may visualize the components of the random vector $(X_1,\ldots,X_k)$ as, approximatively,
normalized random matrices $(M_n(t_1),\ldots,M_n(t_k))$
in (\ref{mn}) with large $n$.

Let $I$ be an arbitrary set. A {\em semicircular family} indexed by $I$ is a collection
of random variables $\{X_i:\,i\in I\}$ such that, for every $k\geq 1$ and every
$(i_1,\ldots,i_k)\in I^k$,
the vector $(X_{i_1},\ldots,X_{i_k})$ has a $k$-dimensional
semicircular distribution.
When $X=\{X_i:\,i\in I\}$ is a semicircular family for which $I$ is denumerable (resp. for which $I=\R_+$),
we say that $X$ is a {\em semicircular sequence} (resp. {\em semicircular process}).

The distribution of any centered semicircular family
$\{X_i:\,i\in I\}$  turns out to be uniquely determined by its covariance function
$\Gamma:I^2\to\R$
given by $\Gamma(i,j)=\varphi(X_iX_j)$.
(This is an easy consequence of \cite[Corollary 9.20]{nicaspeicher}.)
When $I=\mathbb{Z}$, the family is said to be {\em stationary} if $\Gamma(i,j)=\Gamma(|i-j|)$ for all $i,j\in \mathbb{Z}$.

Let $X=\{X_k:\,k\in\Z\}$ be a centered semicircular sequence and consider the linear span $\mathcal{H}$ of $X$,
called the {\em semicircular space associated to} $X$. It is a real separable Hilbert space and, consequently,
there exists an isometry $\Phi:\mathcal{H}\to L^2(\R)$. For any $k\in\Z$, set $e_k = \Phi(X_k)$;
we have, for all $k,l\in\Z$,
\[
\int_\R e_k(x)e_l(x)dx=\varphi(X_kX_l)=\Gamma(k,l).
\]
Thus, since the covariance function $\Gamma$ of $X$ characterizes its distribution, we have
\[
\{X_k:\,k\in\Z\} \,\overset{\rm law}{=}\,\left\{ \int_\R e_k(x)dS(x):\,k\in\Z\right\},
\]
with the notation of Section \ref{wigner1st}.

\subsection{Multiple Wigner integral} \label{wigner}

Let $S=\{S(t)\}_{t\in\R}$ be a two-sided free Brownian motion, and let $p\geq 1$ be an integer.
When $f$ belongs to $L^2(\R^p)$ (recall from Section \ref{wigner1st} that it means, in particular, that $f$ is real-valued), we write $f^*$ to indicate the function of $L^2(\R^p)$ given by
$f^*(t_1,\ldots,t_p)=f(t_p,\ldots,t_1)$.

Following \cite{bianespeicher}, let us quickly sketch out the construction of the {\em multiple Wigner
integral} of $f$ with respect to $S$.
Let $\Delta^q\subset\R^q$ be the collection of all diagonals, i.e.
\begin{equation}\label{dp}
\Delta^q=\{(t_1,\ldots,t_q)\in\R^q:\,t_i=t_j\mbox{ for some $i\neq j$}\}.
\end{equation}
For a characteristic function $f={\bf 1}_A$, where $A\subset\R^q$ has the form
$
A=[u_1,v_1]\times\ldots\times [u_q,v_q]
$
with $A\cap \Delta^q=\emptyset$, the $q$th multiple Wigner integral of $f$
is defined
by
\[
I^{S}_q(f)=(S(v_1)-S(u_1))\ldots (S(v_q)-S(u_q)).
\]
We then extend this definition by linearity to simple functions of the form
$
f=\sum_{i=1}^k \alpha_i {\bf 1}_{A_i},
$
where
$
A_i=[u^i_1,v^i_1]\times\ldots\times [u^i_q,v^i_q]
$
are disjoint $q$-dimensional rectangles as above which do not meet the diagonals.
Simple computations show that
\begin{eqnarray}
\ff(I^{S}_q(f))&=&0\label{isomfree1}\\
\ff(I^{S}_q(f)I^{S}_q(g))&=&\langle f,g^*\rangle_{L^2(\R^q)}.\label{isomfree}
\end{eqnarray}
By approximation, the definition of $I^{S}_q(f)$ is extended to all $f\in L^2(\R^q)$,
and (\ref{isomfree1})-(\ref{isomfree}) continue to hold in this more general setting.
If one wants $I^{S}_q(f)$ to be a random variable in the sense of Section \ref{s:rv},
it is necessary that $f$ be {\em mirror symmetric}, that is, $f=f^*$, in order
to ensure that $I^{S}_q(f)$ is self-adjoint, namely $(I^{S}_q(f))^*=I^{S}_q(f)$ (see \cite{knps}).
Observe that $I^{S}_1(f)=\int_\R f(x)dS(x)$ (see Section \ref{wigner1st}) when $q=1$.
We have moreover
\begin{equation}\label{isom-dif-free}
\ff(I^{S}_p(f)I^{S}_q(g))=0\,\,\mbox{ when $p\neq q$, $f\in L^2(\R^p)$ and $g\in L^2(\R^q)$}.
\end{equation}

When $r\in\{1,\ldots,p\wedge q\}$, $f\in L^2(\R^p)$ and
$g\in L^2(\R^q)$,
let us  write $f\overset{r}\frown g$ to indicate the $r$th {\em contraction} of $f$ and $g$, defined as
being the element of $L^2(\R^{p+q-2r})$ given by
\begin{eqnarray}
&&f\overset{r}\frown g (t_1,\ldots,t_{p+q-2r})\label{contr}\\
&=&\int_{\R^r}f(t_1,\ldots,t_{p-r},x_1,\ldots,x_r)g(x_r,\ldots,x_1,t_{p-r+1},\ldots,t_{p+q-2r})dx_1\ldots dx_r.\notag
\end{eqnarray}
By convention, set $f\overset{0}{\frown} g= f\otimes g$ as
being the tensor product of $f$ and $g$.
Since $f$ and $g$ are not necessarily symmetric functions, the position
of the identified variables $x_1,\ldots,x_r$ in (\ref{contr}) is important, in contrast
to what happens in classical probability (see \cite[Section 6.2]{peccati-taqqu:2011-book}).
Observe moreover that
\begin{equation}\label{cauchy}
\|f\overset{r}{\frown} g\|_{L^2(\R^{p+q-2r})}\leq \|f\|_{L^2(\R^p)}\|g\|_{L^2(\R^q)}
\end{equation}
by Cauchy-Schwarz, and also that
$f\overset{p}{\frown} g=\langle f,g^*\rangle_{L^2(\R^p)}$ when $p=q$.

We have the following {\em product formula} (see \cite[Proposition 5.3.3]{bianespeicher}), valid for any $f\in L^2(\R^p)$ and $g\in L^2(\R^q)$:
\begin{equation}\label{prodfree}
I^{S}_p(f)I^{S}_q(g)=\sum_{r=0}^{p\wedge q} I^{S}_{p+q-2r}(f\overset{r}{\frown}g).
\end{equation}
We deduce (by a straightforward induction) that, for any $e\in L^2(\R)$ and any $q\geq 1$,
\begin{equation}\label{lien-tch}
U_q\left(\int_\R e(x)dS_x\right) = I^{S}_q(e^{\otimes q}),
\end{equation}
where $U_0(x)=1$, $U_1(x)=x$, $U_2(x)=x^2-1$, $U_3(x)=x^3-2x$, $\ldots$, is the
sequence of Tchebycheff polynomials of second kind
(determined by the recursion $xU_k=U_{k+1}+U_{k-1}$),
$\int_\R e(x)dS(x)$ is understood as a Wigner integral (as defined in Section \ref{wigner1st}),
and $e^{\otimes q}$ is the $q$th tensor product of $e$, i.e., the symmetric element
of $L^2(\R^q)$ given by $e^{\otimes q}(t_1,\ldots,t_q)=e(t_1)\ldots e(t_q)$.

\section{Non-commutative Tchebycheff processes and basic properties}\label{s:ncTp}

We define first the non-commutative fractional Brownian motion, then the non-commutative Rosenblatt process,
which is a multiple Wigner integral of order 2, and then we introduce the general non-commutative Tchebycheff processes
involving multiple Wigner integrals of arbitrary order.

\subsection{Non-commutative fractional Brownian motion}

The classical fractional Brownian motion was introduced by Kolmogorov \cite{kolmo} and developed by Mandelbrot and
Van Ness \cite{mandel}.

\begin{defi}\label{d:ncfbm}
Let $H\in(0,1)$. A non-commutative fractional Brownian motion (ncfBm in short) of Hurst parameter $H$ is a centered semicircular process
$S_H=\{S_H(t):\,t\geq 0\}$
with covariance function
\begin{equation}\label{e:cov}
\varphi(S_H(t)S_H(s))=\frac12\big(t^{2H}+s^{2H}-|t-s|^{2H}\big).
\end{equation}
\end{defi}
It is readily checked that $S_{1/2}$ is nothing but a one-sided free Brownian motion. Immediate properties of $S_H$, proved
in Corollary \ref{c:ss},
include the selfsimilarity property (that is, for all $a>0$ the process $\{a^{-H}S_H(at):\,t\geq 0\}$ is a ncfBm of parameter $H$)
and the stationary property of the increments (that is, for all $h>0$ the process $\{S_H(t+h)-S_H(h):\,t\geq 0\}$ is a ncfBm of parameter $H$).
Conversely, ncfBm of parameter $H$ is the only standardized semicircular process to verify these two properties, since
they determine the covariance (\ref{e:cov}).

It is interesting to notice that ncfBm may be easily represented as a Wigner integral as follows:
\[
S_H\overset{\rm f.d.d.}{=}\sqrt{\frac{2H}{(H-\frac12)\beta(H-\frac12,2-2H)}}
\int_\R \left( (\cdot-x)_+^{H-\frac12}-(-x)_+^{H-\frac12}\right)dS(x).
\]
Here, $\beta$ stands for the usual Beta function.
In the classical probability case, one has a similar representation with $S$ replaced
by a Brownian motion.

As an illustration, we will now show that normalized sums of semicircular sequences can converge to ncfBm.
Let $\{X_k:\,k\in\Z\}$ be a stationary semicircular sequence with $\ff(X_k)=0$ and $\ff(X_k^2)=1$,
and suppose that its correlation kernel $\rho(k-l)=\ff(X_kX_l)$ verifies
\begin{equation}\label{rho2}
\sum_{k,l=1}^n\rho(k-l)\sim Kn^{2H}L(n)\quad\mbox{as $n\to\infty$},
\end{equation}
with $L:(0,\infty)\to (0,\infty)$ slowly varying at infinity, $0<H<1$ and $K$ a positive constant.
Consider the non-commutative stochastic process
\[
Z_n(t)=
\frac1{n^{H}\sqrt{L(n)}}
\sum_{k=1}^{[nt]}X_k,\quad t\geq 0.
\]
For any $t\geq s\geq 0$, we have, as $n\to\infty$,
\begin{eqnarray*}
&&\ff\left[Z_n(t)Z_n(s)\right]\\
&=& \frac12\ff\left[Z_n(t)^2\right] +
\frac12\ff\left[Z_n(s)^2\right]
-
\frac12 \ff\left[(Z_n(t)-Z_n(s))^2\right]\\
&=&\frac{1}{2n^{2H}L(n)}\sum_{i,j=1}^{[nt]}\rho(i-j) +
\frac{1}{2n^{2H}L(n)}\sum_{i,j=1}^{[ns]}\rho(i-j)
-\frac{1}{2n^{2H}L(n)}\sum_{i,j=[ns]+1}^{[nt]}\rho(i-j)\\
&=&\frac{1}{2n^{2H}L(n)}\sum_{i,j=1}^{[nt]}\rho(i-j) +
\frac{1}{2n^{2H}L(n)}\sum_{i,j=1}^{[ns]}\rho(i-j)
-\frac{1}{2n^{2H}L(n)}\sum_{i,j=1}^{[nt]-[ns]}\rho(i-j)\\
&\to&\frac{K}2\big(t^{2H}+s^{2H}-(t-s)^{2H}\big)=K\,\varphi(S_H(t)S_H(s)).
\end{eqnarray*}
Let $p\geq 1$ as well as $t_1,\ldots,t_p\geq 0$.
Since the $X_k$'s are centered and jointly semicircular, the process $Z_n$ is centered and semicircular as well,
and we have shown that, as $n\to\infty$,
\[
Z_n\overset{{\rm f.d.d.}}{\to} \sqrt{K}\,S_{H}.
\]

\subsection{Non-commutative Rosenblatt process}

The Hermite process indexed by $q\geq 1$ appeared as a limit in Theorem \ref{main-cl} in the classical probability setting.
When $q=1$, it is fractional Brownian motion. When $q=2$, it is the Rosenblatt process, introduced in Taqqu \cite{Taqqu1} and
which appears as a limit in many statistical tests.
See Taqqu \cite{taqqu:2011} for a recent overview.
We introduce here the non-commutative Rosenblatt process.

\begin{defi}\label{d:rosen}
Let $H\in(\frac12,1)$. The non-commutative Rosenblatt process of parameter $H$ is the non-commutative stochastic process
defined by the double Wigner integral
\begin{equation}\label{e:rosen}
R_H(t)=R_{H,2}(t)=I^{S}_2\big(f_{H}(t,\cdot)\big),\quad t\geq 0,
\end{equation}
where
\begin{equation}\label{e:kernelrosen}
f_{H}(t,x,y)=
\frac{\sqrt{H(2H-1)}}{
\beta(\frac{H}2,1-H)
}
\int_0^t (s-x)^{\frac{H}2-1}_+(s-y)^{\frac{H}2-1}_+ds,
\end{equation}
with $\beta$ the usual Beta function.
\end{defi}
Using the relations
\[
\int_{\R} (t-x)^{\frac{H}2-1}_+(s-x)^{\frac{H}2-1}_+dx = \beta(\frac{H}2,1-H)|t-s|^{H-1},
\]
as well as
\begin{equation}\label{relevant}
H(2H-1)\iint_{[0,T]^2}|t-s|^{2H-2} dsdt= T^{2H},\quad T>0,
\end{equation}
it is straightforward to check that
$\varphi(R_{H}(t)^2)=\int_{\R^2}f_{H}(t,x,y)^2dxdy=t^{2H}$.

The Rosenblatt process at  (fixed)  time $t$ is a double Wigner integral whose kernel $f_H(t,\cdot)$
is symmetric, see (\ref{e:kernelrosen}). As such, it enjoys useful properties, that we derive now in full generality.
Assume then that $f\in L^2(\R^2)$ is a given symmetric kernel.
One of the most effective ways of dealing with $I^{S}_2(f)$
is to associate to $f$ the following Hilbert-Schmidt operator:
\begin{equation}
A_{f}:L^2(\R)\to L^2(\R);\,\,g\mapsto \int_\R f(\cdot,y)g(y)dy.
\label{Hilsch}
\end{equation}%
In other words, $A_{f}$ transforms an element $g$ of $L^2(\R)$ into the contraction $f\overset{1}\frown g\in
L^2(\R)$. We write
$\left\{ \lambda _{f,j} : j\geq 1 \right\}$ and $\left\{
e_{f,j} : j\geq 1\right\}$, respectively, to indicate the
eigenvalues of $A_{f}$ and the corresponding eigenvectors (forming
an orthonormal system in $L^2(\R)$).

Some useful relations between all these objects
are explained in the next proposition. The proof, which is omitted here, relies on elementary functional analysis
(see e.g. Section 6.2 in \cite{HirschLacombe}).

\begin{prop}
Let $f$ be a symmetric element of $L^2(\R^2)$, and let the above notation prevail.
\begin{enumerate}
\item The series $\sum_{j=1}^{\infty }\lambda _{f,j}^{p}$
converges for every $p\geq2$, and $f$ admits the expansion
\begin{equation}
f=\sum_{j=1}^{\infty }\lambda _{f,j} \,\,\,\big(e_{f,j}\otimes e_{f,j}\big)\text{,} \label{squ}
\end{equation}%
where the convergence takes place in $L^2(\R^2)$.
\item For every $p\geq 2$, one has the relations%
\begin{equation}
{\rm Tr} ( A_{f}^{p}) =
\int_{\R^{p}}f(x_1,x_2)\ldots f(x_{p-1},x_p)f(x_p,x_1)dx_1\ldots dx_{p}
=\sum_{j=1}^{\infty
}\lambda _{f,j}^{p}\text{,} \label{capz}
\end{equation}%
where ${\rm Tr} ( A_{f}^{p})$ stands for the trace of the $p$th
power of $A_{f}$.
\end{enumerate}
\end{prop}

\bigskip

In the following statement we collect some  facts
concerning the law of a random variable of the type
$I^{S}_{2}\left( f\right)$.

\begin{prop}\label{P : 2cum}
Let $F=I^{S}_{2}\left( f\right) $, where $f$ is a symmetric element of $L^2(\R^2)$.
\begin{enumerate}

\item The following equality holds:%
\begin{equation}
F\,\,\overset{\rm Law}{=}\,\,\sum_{j=1}^{\infty }\lambda _{f,j}\left( S_{j}^{2}-1\right)
\text{,} \label{expX2}
\end{equation}%
where $( S_{j})_{j\geq 1}$ is a sequence of
freely independent $\mathcal{S}(0,1)$ random variables, and the series
converges in $L^2(\mathscr{A},\varphi)$.

\item For every $p\geq 2$, the $p$th free cumulant of $F$ is given by the following formula:
\begin{equation}
\kappa _{p}\left( F,\ldots,F\right)
=\int_{\R^{p}}f(x_1,x_2)\ldots f(x_{p-1},x_p)f(x_p,x_1)dx_1\ldots dx_{p}
\text{.} \label{cumX2}
\end{equation}
\end{enumerate}
\end{prop}

\noindent\textit{Proof. } Relation
(\ref{expX2}) is an immediate consequence of (\ref{squ}), of the identity \[%
I^{S}_{2}\left( e_{f,j}\otimes e_{f,j}\right) = I^{S}_{1}\left( e_{f,j}\right)
^{2}-1, \] as well as of the fact that the $\left\{ e_{f,j}\right\} $ are
orthonormal (implying that the sequence $\{I^{S}_{1}\left(
e_{f,j}\right) : j\geq 1 \}$ is composed of freely independent $\mathcal{S}(0,1)$ random variables). To
prove (\ref{cumX2}),
it suffices to use the linearization property (\ref{cum}) of free cumulants, as well as the fact\footnote{Since all the free cumulants of $S^2$ are equal to 1, $S^2$ has a free Poisson law with mean 1.}
 that
\[
\kappa_p(S^2-1,\ldots,S^2-1)=\kappa_p(S^2,\ldots,S^2)+\kappa_p(-1,\ldots,-1)=\kappa_p(S^2,\ldots,S^2)=1
\]
for all $p\geq 2$, see e.g. \cite[Proposition 12.13]{nicaspeicher} for the last equality.
We thus obtain the desired conclusion by means of (\ref{capz}).
\qed

\bigskip

In the classical probability setting where $S_j$ is $\mathcal{N}(0,1)$, there is an additional
factor of $2^{p-1}(p-1)!$ in (\ref{cumX2}). See e.g. Taqqu \cite{Taqqu1}.

\subsection{Non-commutative Tchebycheff processes}

The classical probability versions of these processes are the Hermite processes.
See e.g. Peccati-Taqqu \cite[Section 9.5]{peccati-taqqu:2011-book}.

\begin{defi}\label{d:tcheb}
Let $H\in(\frac12,1)$. The $q$th non-commutative Tchebycheff process of parameter $H$ is the non-commutative stochastic process defined by the Wigner integral
\begin{equation}\label{e:tcheb}
R_{H,q}(t)=I^{S}_q\big(f_{H,q}(t,\cdot)\big),\quad t\geq 0,
\end{equation}
where
\begin{equation}\label{e:kernel}
f_{H,q}(t,x_1,\ldots,x_q)=
\frac{\sqrt{H(2H-1)}}{
\beta(\frac12-\frac{1-H}q,\frac{2-2H}q)^{q/2}
}
\int_0^t (s-x_1)^{-(\frac12+\frac{1-H}q)}_+\ldots(s-x_q)^{-(\frac12+\frac{1-H}q)}_+ds,
\end{equation}
with $\beta$ the usual Beta function.
\end{defi}

The process $R_{H,q}$ becomes the non-commutative fractional Brownian motion and the non-commutative
Rosenblatt process when $q=1$ and $q=2$ respectively. Note however that when $q=1$ the process is defined
for $H$ between 0 and 1.

Using the relations
\begin{equation}\label{gamma}
\int_{\R} (t-x)^{\gamma-1}_+(s-x)^{\gamma-1}_+dx = \beta(\gamma,1-2\gamma)|t-s|^{2\gamma-1},\quad 0<\gamma<1/2,
\end{equation}
and (\ref{relevant}),
we easily get that
\begin{equation}\label{335}
\varphi(R_{H,q}(t)^2)=\int_{\R^q}f_{H,q}(t,x_1,\ldots,x_q)^2dx_1\ldots dx_q=t^{2H}.
\end{equation}

The process $R_{H,q}$ has stationary increments and is selfsimilar with parameter $H$, as stated
in Corollary \ref{c:ss} below.

\section{Computation of joint moments}\label{s:joint}

Let $p\geq 2$ be a given integer.
Let $f_1,\ldots,f_p$ be real functions of $q_1,\ldots,q_p$ variables respectively.
Write ${\bf q}_p=(q_1,\ldots,q_p)$.
We want to compute
\begin{equation}\label{joint}
(\ldots((f_1\overset{r_1}\frown f_2)\overset{r_2}\frown f_3)\ldots)\overset{r_{p-1}}\frown f_p
\end{equation}
for some functions $f_1,\ldots,f_p$ of interest.
The contraction operator $\overset{r}\frown$ is defined in (\ref{contr}).
The expression (\ref{joint}) makes sense if and only if ${\bf r}=(r_1,\ldots,r_{p-1})\in A({\bf q}_p)$, where
$A({\bf q}_p)$ is the set of those $(r_1,\ldots,r_{p-1})\in \{0,\ldots,q_2\}
\times\ldots\times
\{0,\ldots,q_p\}$
such that
\begin{eqnarray}
&&r_1\leq q_1,\,\,r_2\leq q_1+q_2-2r_1,\,\,r_3\leq q_1+q_2+q_3-2r_1-2r_2,\notag\\
&&\hskip3cm \ldots,\,\,r_{p-1}\leq q_1+\ldots+q_{p-1}-2r_1-\ldots-2r_{p-2},
\label{ap}
\end{eqnarray}
and (\ref{joint}) equals a real number if and only if ${\bf r}\in B({\bf q}_p)$, where
\begin{equation}\label{bp}
B({\bf q}_p)=\big\{{\bf r}=(r_1,\ldots,r_{p-1})\in A({\bf q}_p):\,2r_1+\ldots+2r_{p-1}=q_1+\ldots+q_p\big\}.
\end{equation}

Indeed, for $f_1\overset{r_1}\frown f_2$ to make sense, we need $0\leq r_1\leq q_1\wedge q_2$. Then,
for $(f_1\overset{r_1}\frown f_2)\overset{r_2}\frown f_3$ to make sense, we need
$0\leq r_2\leq q_3\wedge(q_1+q_2-2r_1)$; this is because $f_1\overset{r_1}\frown f_2$ has $q_1+q_2-2r_1$
variables and $f_3$ has $q_3$ variables. It follows, by induction, that for (\ref{joint})
to make sense, it is necessary that $(r_1,\ldots,r_{p-1})\in A({\bf q}_p)$.
In order for (\ref{joint}) to be a scalar, we need $(r_1,\ldots,r_{p-1})\in B({\bf q}_p)$, since the number
of variables of (\ref{joint}) is given by $q_1+\ldots+q_p-2r_1-\ldots-2r_{p-1}$.

The following lemma gives the value of (\ref{joint}) for functions $f_1,\ldots,f_p$ of interest.
The result involves an array of non-negative integers
\begin{equation}\label{alphaij}
\alpha_{ij}({\bf r}),\quad 1\leq i<j\leq p,
\end{equation}
which are defined as follows for  ${\bf r}\in B({\bf q}_p)$.
Consider the following figure, where there are $q_1+\ldots+q_p$ dots.
The first $q_1$ corresponds
to the $q_1$ variables of $f_1$, $\ldots$,  the last $q_p$ dots corresponds to the $q_p$ variables
of $f_p$.

\bigskip
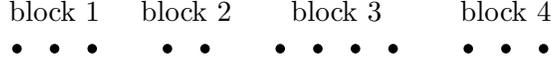
\begin{figure}[!ht]
\begin{center}
\begin{pspicture}(-3,0)(4,.5)
\psdot(-3,0)
\psdot(-2.5,0)
\psdot(-2,0)
\psdot(-1,0)
\psdot(-0.5,0)
\psdot(.5,0)
\psdot(1,0)
\psdot(1.5,0)
\psdot(2,0)
\psdot(3,0)
\psdot(3.5,0)
\psdot(4,0)
\rput(-2.5,.5){block $1$}
\rput(-0.75,.5){block $2$}
\rput(1.25,.5){block $3$}
\rput(3.5,.5){block $4$}
\end{pspicture}
\label{fig1}
\caption{$p=4$, $q_1=3$, $q_2=2$, $q_3=4$ and $q_4=3$}
\end{center}
\end{figure}

We shall associate pairs of dots according to the mechanism described below, and we say
that a dot is {\em available} if it has not been associated so far. The association rule, call it $(A)$,
involves associating the left most available dot in a block with the right most available dot in
preceding blocks.

To perform these associations, proceed as follows. Start with block $j=2$ and
do $r_1$ associations with block $1$, following the association rule $(A)$.

\bigskip
\begin{figure}[!ht]
\begin{center}
\begin{pspicture}(-3,-.25)(4,.5)
\psdot(-3,0)
\psdot(-2.5,0)
\psdot(-2,0)
\psdot(-1,0)
\psdot(-0.5,0)
\psdot(.5,0)
\psdot(1,0)
\psdot(1.5,0)
\psdot(2,0)
\psdot(3,0)
\psdot(3.5,0)
\psdot(4,0)
\psline(-2,0)(-2,-0.25)
\psline(-2,-.25)(-1,-.25)
\psline(-1,0)(-1,-0.25)
\rput(-2.5,.5){block $1$}
\rput(-0.75,.5){block $2$}
\rput(1.25,.5){block $3$}
\rput(3.5,.5){block $4$}
\end{pspicture}
\label{fig2}
\caption{$r_1=1$}
\end{center}
\end{figure}

Proceed to block $j=3,4,\ldots,p$. In block $j\geq 3$, associate $r_{j-1}$ dots with available dots
in the preceding blocks following the association rule $(A)$. Once block $j=p$ is done, all
dots have been associated pairwise.

Figure 3 below illustrates an example of associations.

\bigskip
\begin{figure}[!ht]
\begin{center}
\begin{pspicture}(-3,-.75)(4,.5)
\psdot(-3,0)
\psdot(-2.5,0)
\psdot(-2,0)
\psdot(-1,0)
\psdot(-0.5,0)
\psdot(.5,0)
\psdot(1,0)
\psdot(1.5,0)
\psdot(2,0)
\psdot(3,0)
\psdot(3.5,0)
\psdot(4,0)
\psline(-2,0)(-2,-0.25)
\psline(-2,-.25)(-1,-.25)
\psline(-1,0)(-1,-0.25)
\psline(2,0)(2,-0.25)
\psline(2,-.25)(3,-.25)
\psline(3,0)(3,-0.25)
\psline(.5,0)(.5,-0.25)
\psline(.5,-.25)(-.5,-.25)
\psline(-.5,0)(-.5,-0.25)
\psline(1,0)(1,-0.5)
\psline(-2.5,-.5)(1,-.5)
\psline(-2.5,0)(-2.5,-0.5)
\psline(3.5,0)(3.5,-0.5)
\psline(1.5,-.5)(3.5,-.5)
\psline(1.5,0)(1.5,-0.5)
\psline(4,0)(4,-0.75)
\psline(4,-.75)(-3,-.75)
\psline(-3,0)(-3,-0.75)
\rput(-2.5,.5){block $1$}
\rput(-0.75,.5){block $2$}
\rput(1.25,.5){block $3$}
\rput(3.5,.5){block $4$}
\end{pspicture}
\label{fig3}
\caption{$p=4$, $q_1=3$, $q_2=2$, $q_3=4$, $q_4=3$, $r_1=1$, $r_2=2$ and $r_3=3$}
\end{center}
\end{figure}
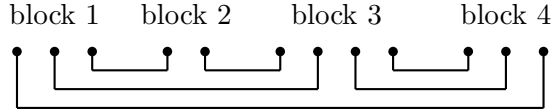

\begin{defi}\label{def:alpha}
When ${\bf r}=(r_1,\ldots,r_{p-1})\in B({\bf q}_p)$,
denote by $\alpha_{ij}({\bf r})$ the number of associations between dots of block $j$ with dots of block $i$, $1\leq i<j\leq p$.
\end{defi}
For instance, in Figure 3, we have
\[
\left(
\begin{array}{ccc}
\alpha_{12}({\bf r})&\alpha_{13}({\bf r})&\alpha_{14}({\bf r})\\
&\alpha_{23}({\bf r})&\alpha_{24}({\bf r})\\
&&\alpha_{34}({\bf r})
\end{array}
\right)
=\left(
\begin{array}{ccc}
1&1&1\\
&1&0\\
&&2
\end{array}
\right)
.
\]

The following lemma gives an explicit expression for (\ref{joint}) for functions $f_1,\ldots,f_q$ which
appear in the sequel.

\begin{lemma}\label{quivaienb}
Let $T\subset\R$, and let $e:T\times \R\to \R$ be a measurable function.
Fix also integers $p\geq 2$ and $q_1,\ldots,q_p\geq 1$, and let $\nu_1,\ldots,\nu_p$ be given signed measures on $\R$.
Assume further that
\[
\int_{\R^q} \left(\int_T |e(s,x_1)|\ldots |e(s,x_{q_i})|\,|\nu_i|(ds)\right)^2dx_1\ldots dx_{q_i}<\infty, \quad i=1,\ldots,p.
\]
Finally, for any $i=1,\ldots,p$,
define $f_i:\R^{q_i}\to \R$ to be
\[
f_i(x_1,\ldots,x_{q_i})=\int_T e(s,x_1)\ldots e(s,x_{q_i})\nu_i(ds).
\]
Then, for any ${\bf r}=(r_1,\ldots,r_{p-1})\in B({\bf q}_p)$ (recall the definition (\ref{bp}) of $B({\bf q}_p)$), we have
\[
(\ldots((f_1\overset{r_1}\frown f_2)\overset{r_2}\frown f_3)\ldots)\overset{r_{p-1}}\frown f_p
=\int_{T^p}
\prod_{1\leq i<j\leq p} \left(
\int_\R e(s_i,x)e(s_j,x)dx
\right)^{\alpha_{ij}({\bf r})}
\!\!\!\!\!
 \nu_1(ds_1)\ldots \nu_p(ds_p).
\]
\end{lemma}
\noindent
{\it Proof}.
Fix  ${\bf r}=(r_1,\ldots,r_{p-1})\in B({\bf q}_p)$.
We have
\begin{eqnarray*}
&&
(\ldots((f_1\overset{r_1}\frown f_2)\overset{r_2}\frown f_3)\ldots)\overset{r_{p-1}}\frown f_p
\\
&=&\int_{T^p}
(\ldots((e(s_1,\cdot)^{\otimes q_1}\overset{r_1}\frown e(s_2,\cdot)^{\otimes q_2})\overset{r_2}\frown
e(s_3,\cdot)^{\otimes q_3})\ldots)\overset{r_{p-1}}\frown
e(s_p,\cdot)^{\otimes q_p}
\nu_1(ds_1)\ldots \nu_p(ds_p),
\end{eqnarray*}
where we expressed $f_i$ as $f_i(\cdot)=\int_T e(s,\cdot)^{\otimes q_i}\nu_i(ds)$.
What matters in the computation of
$(\ldots((e(s_1,\cdot)^{\otimes q_1}\overset{r_1}\frown e(s_2,\cdot)^{\otimes q_2})\overset{r_2}\frown
e(s_3,\cdot)^{\otimes q_3})\ldots)\overset{r_{p-1}}\frown
e(s_p,\cdot)^{\otimes q_p}$
is, for each $i<j$, the number $\alpha_{ij}({\bf r})$ of associations between the block $j$ and the block $i$.
Hence
\begin{eqnarray*}
&&(\ldots((e(s_1,\cdot)^{\otimes q_1}\overset{r_1}\frown e(s_2,\cdot)^{\otimes q_2})\overset{r_2}\frown e(s_3,\cdot)^{\otimes q_3})\ldots)\overset{r_{p-1}}\frown
e(s_p,\cdot)^{\otimes q_p}\\
&&\hskip6cm=\prod_{1\leq i<j\leq p} \left(
\int_\R e(s_i,x)e(s_j,x)dx
\right)^{\alpha_{ij}({\bf r})},
\end{eqnarray*}
and the desired formula follows.
\qed

\medskip

We shall apply Lemma \ref{quivaienb} with a discrete $\nu=\sum \delta_k$ where $\delta_k$ is Dirac mass, and with a
Lebesgue-type $\nu$.
The following result plays an important role
in the proof of the non-central limit theorem (\ref{tcl2}).

\begin{prop}\label{l:key}
Fix two integers $q\geq 1$ and $p\geq 2$, and let $t_1,\ldots,t_p$ be positive real numbers.
Set ${\bf q}_p=(q_1,\ldots,q_p)=(q,\ldots,q)$.
Recall the definition (\ref{alphaij}) of $\alpha_{ij}({\bf r})$ when ${\bf r}=(r_1,\ldots,r_{p-1})\in B({\bf q}_p)$.
Then, the following two assertions hold:

1. Let $X=\{X_k:\,k\in\Z\}$ be a stationary semicircular sequence with $\ff(X_k)=0$ and $\ff(X_k^2)=1$,
and let $\rho(k-l)=\ff(X_kX_l)$ be its correlation kernel.
Write $U_q$ to indicate the $q$th Tchebycheff polynomial, and define
\begin{equation}\label{gnn}
G_n(t)=\sum_{k=1}^{[nt]}U_q(X_k), \quad t\geq 0.
\end{equation}
Then
\begin{equation}\label{kthmomentdiscret}
\varphi(G_n(t_1)\ldots G_n(t_p))=\sum_{k_1=1}^{[nt_1]}\ldots \sum_{k_p=1}^{[nt_p]}
\,\,\,\sum_{(r_1,\ldots,r_{p-1})\in B({\bf q}_p)} \,\,\, \prod_{1\leq i<j\leq p}\rho(k_i-k_j)^{\alpha_{ij}({\bf r})}.
\end{equation}

2. Let $R_{H,q}$ be a $q$th non-commutative Tchebycheff process of parameter $H\in(1/2,1)$.
Then
\begin{eqnarray}\label{kthmoment2}
&&\varphi(R_{H,q}(t_1)\ldots R_{H,q}(t_p))\\
&=&
H^{p/2}(2H-1)^{p/2}
\int_0^{t_1}ds_1\ldots\int_0^{t_p}ds_p\,\,\,\sum_{(r_1,\ldots,r_{p-1})\in B({\bf q}_p)} \,\,\,
\prod_{1\leq i<j\leq p}|s_i-s_j|^{-\alpha_{ij}({\bf r})\times\frac{2-2H}{q}}.
\notag
\end{eqnarray}
\end{prop}
\noindent
{\it Proof}.
1. Let $\{e_k:\,k\in\Z\}\subset L^2(\R)$ be defined as in Section \ref{sec:semfree}.
In the definition (\ref{gnn}) of $G_n$, we can assume without loss of generality that
\[
X_k=I^{S}_1(e_k)=\int_\R e_k(x)dS(x), \quad k\in\Z.
\]
By virtue of (\ref{lien-tch}), we have $G_n(t)=I^{S}_q(g_n(t,\cdot))$,
with
\[
g_n(t,x_1,\ldots,x_q) = \sum_{k=1}^{[nt]} e_k(x_1)\ldots e_k(x_q).
\]
By the product formula (\ref{prodfree}),
\[
G_n(t_1)G_n(t_2)=I^{S}_q\big(g_n(t_1,\cdot)\big)I^{S}_q\big(g_n(t_2,\cdot)\big)
=\sum_{r_1=0}^q I^{S}_{2q-2r_1}\big(g_n(t_1,\cdot)\overset{r_1}{\frown}g_n(t_2,\cdot)\big).
\]
Iterative applications of the product formula (\ref{prodfree}) lead to
\begin{eqnarray}\label{kthmomentwithoutexpectation}
&&G_n(t_1)\ldots G_n(t_p)\\
&=&\!\!\!\!\sum_{(r_1,\ldots,r_{p-1})\in A({\bf q}_p)}\!\!\!\!
I^{S}_{pq-2r_1-\ldots-2r_{p-1}}\big(
(\ldots((g_n(t_1,\cdot)\overset{r_1}\frown g_n(t_2,\cdot))\overset{r_2}\frown g_n(t_3,\cdot))\ldots)
\overset{r_{p-1}}\frown g_n(t_p,\cdot)
\big),\notag
\end{eqnarray}
where $A({\bf q}_p)$ is defined by (\ref{ap}).
By applying $\varphi$ to (\ref{kthmomentwithoutexpectation}) and using (\ref{isomfree1}), we deduce that
the only non-zero terms occur when $pq-2r_1-\ldots-2r_{p-1}=0$. Since $I_0(c)=c$, $c\in\R$, we get
\begin{eqnarray}\label{kthmoment}
&&\varphi(G_n(t_1)\ldots G_n(t_p))\\
&=&\sum_{(r_1,\ldots,r_{p-1})\in B({\bf q}_p)}
(\ldots((g_n(t_1,\cdot)\overset{r_1}\frown g_n(t_2,\cdot))\overset{r_2}\frown g_n(t_3,\cdot))\ldots)
\overset{r_{p-1}}\frown g_n(t_p,\cdot)
,\notag
\end{eqnarray}
with $B({\bf q}_p)$ defined by (\ref{bp}).
Now, recall from Section \ref{sec:semfree} the following property of $e_k$:
\begin{equation}\label{identityquivabien}
\int_\R e_k(x)e_l(x)dx = \varphi(X_kX_l)=\rho(k-l).
\end{equation}
Hence, using Lemma \ref{quivaienb}
with $T=\Z$,
$\nu_i=\sum_{k=1}^{[nt_i]}\delta_k$ ($\delta_k$ being the Dirac mass at $k$)
and $e(k,x)=e_k(x)$,
we obtain the formula (\ref{kthmomentdiscret})
for $\varphi(G_n(t_1)\ldots G_n(t_p))$.\\

2. Reasoning as in Point 1 above, we get here that
\begin{eqnarray}
&&\varphi(R_{H,q}(t_1)\ldots R_{H,q}(t_p))\label{star2}\\
&=&\sum_{(r_1,\ldots,r_{p-1})\in B({\bf q}_p)}
(\ldots((f_{H,q}(t_1,\cdot)\overset{r_1}\frown f_{H,q}(t_2,\cdot))\overset{r_2}\frown f_{H,q}(t_3,\cdot))\ldots)
\overset{r_{p-1}}\frown f_{H,q}(t_p,\cdot)
,\notag
\end{eqnarray}
with $B({\bf q}_p)$ defined by (\ref{bp}).
Lemma \ref{quivaienb},
with $T=\R_+$, measures
\[\nu_i(ds)=\frac{\sqrt{H(2H-1)}}{\beta(\frac12-\frac{1-H}q,\frac{2-2H}q)^{q/2}}{\bf 1}_{[0,t_i]}(s)ds
\]
and the function $e(s,x)=(s-x)_+^{-(\frac12+\frac{1-H}q)}$,
then yields, because of (\ref{gamma}),
\begin{eqnarray}\label{kthmomentfinal2}
&&(\ldots((f_{H,q}(t_1,\cdot)\overset{r_1}\frown f_{H,q}(t_2,\cdot))\overset{r_2}\frown f_{H,q}(t_3,\cdot))\ldots)
\overset{r_{p-1}}\frown f_{H,q}(t_p,\cdot)\\
&=&H^{p/2}(2H-1)^{p/2}\,\int_0^{t_1}ds_1\ldots\int_0^{t_p}ds_p\,\,\,\prod_{1\leq i<j\leq p}
|s_i-s_j|^{-\alpha_{ij}({\bf r})\times\frac{2-2H}{q}}.\notag
\end{eqnarray}
By inserting (\ref{kthmomentfinal2}) in (\ref{star2}), we obtain the
formula (\ref{kthmoment2}) for $\varphi(R_{H,q}(t_1)\ldots R_{H,q}(t_p))$.\qed

\begin{cor}\label{c:ss}
The non-commutative Tchebycheff process $R_{H,q}$ has stationary increments and is selfsimilar with parameter $H$.
\end{cor}
\noindent
{\it Proof}. Since the law is determined by the moments, it suffices to use expression (\ref{kthmoment2}).
Let $h>0$. Replacing $R_{H,q}(t_i)$ by $R_{H,q}(t_i+h)-R_{H,q}(h)$ in the left-hand side of (\ref{kthmoment2})
changes the integrals $\int_0^{t_i}$ in the right-hand side by integrals $\int_{h}^{t_i+h}$.
Since this does not modify the right-hand side, the process $R_{H,q}$ has stationary increments.
To prove selfsimilarity, let $a>0$, replace each $t_1,\ldots,t_p$ by $at_1,\ldots,at_p$ in (\ref{kthmoment2}) and note that the right-hand side is then
multiplied by a factor $a^{pH}$.
\qed

\section{Proof of Theorem \ref{main}}\label{s:proof}

In the proof of the central limit theorem (\ref{tcl1}), we shall use the following Wiener-Wigner transfer principle, established in
\cite[Theorem 1.6]{NPS}. It provides an equivalence between multidimensional limit theorems involving multiple Wiener integrals and multiple
Wigner integrals respectively, whenever the limits of the multiple Wiener integrals are normal.

\begin{prop}\label{transfer}{\rm (Statement of \cite[Theorem 1.6]{NPS})}
 Let $d\geq 1$ and $q_1,\dots,q_d$ be some fixed integers, and consider a positive definite symmetric matrix
$c =\{c(i,j) : i,j=1,...,d\}$. Let $(G_1,\ldots,G_d)$ be a $d$-dimensional Gaussian vector
and $(\Sigma_1,\dots,\Sigma_d)$ be a semicircular vector, both with covariance $c$
as defined in Section \ref{sec:semfree}.
For each $i=1,\dots,d$, we consider a sequence $\{f_{i,n}\}_{n\geq 1}$ of symmetric functions
in $L^2(\R^{q_i}_+)$.
Let $B$ be a classical Brownian motion and let $I^B_q(\cdot)$ stand for the $q$th multiple Wiener integral.
Let $S$ be a free Brownian motion and let $I^S_q(\cdot)$ stand for the $q$th multiple Wigner integral.
Then:
\begin{enumerate}
\item For all $i,j=1,\ldots,d$ and as $n\to\infty$, $\varphi[I^{S}_{q_i}(f_{i,n})I^{S}_{q_j}(f_{j,n})] \to c(i,j)$
if and only if $E[I^B_{q_i}(f_{i,n})I^B_{q_j}(f_{j,n})] \to \sqrt{(q_i)!(q_j)!}\,c(i,j)$.
\item If the asymptotic relations in {\rm (1)} are verified then, as $n\to\infty$,
\[
\big(I^{S}_{q_1}(f_{1,n}),\dots,I^{S}_{q_d}(f_{d,n})\big)\overset{\rm law}{\to}(\Sigma_1,\dots,\Sigma_d)
\]
if and only if
\[
\big(I^B_{q_1}(f_{1,n}),\dots,I^B_{q_d}(f_{d,n)})\big)\overset{\rm law}{\to}\big(\sqrt{(q_1)!}G_1,\dots,\sqrt{(q_d)!}G_d\big).
\]
\end{enumerate}
\end{prop}

\subsection{Proof of the central limit theorem (\ref{tcl1})}
Recall that $X=\{X_k:\,k\in\Z\}$ is a stationary {\sl semicircular} sequence with $\ff(X_k)=0$, $\ff(X_k^2)=1$ and correlation $\rho$.
Consider first its Gaussian counterpart (in the usual probabilistic sense), namely
$Y=\{Y_k:\,k\in\Z\}$ where $Y$ is a stationary {\sl Gaussian} sequence with mean $0$ and same correlation $\rho$.
When dealing with $Y$, the important polynomials are not the Tchebycheff polynomials but the Hermite polynomials, defined
as $H_0(y)=1,H_1(y)=y,H_2(y)=y^2-1, H_3(y)=y^3-3y,\ldots,$ and determined by the recursion $yH_k=H_{k+1}+kH_{k-1}$.

We assume in this proof that $\sum_{k\in\Z}|\rho(k)|^q<\infty$; this implies
$\sum_{k\in\Z}|\rho(k)|^s<\infty$ for all $s\geq q$.
Since $Q$ given by (\ref{decompo}) is a polynomial, we can choose $N$ large enough so
that $a_s=0$ for all $s\geq N$.
Set
\[
W_n(H_s,t)=\sum_{k=1}^{[nt]}H_s(Y_k), \quad t\geq 0, \quad s=q,\ldots,N.
\]
The celebrated Breuer-Major theorem (see \cite{breuermajor}, see also \cite[Chapter 7]{ivangiobook} for a
modern proof, and see Theorem \ref{main-cl}, part 1, for the statement) asserts that
\begin{equation}\label{h_nu}
\left(
\frac{W_n(H_q,\cdot)}{\sqrt{n}},\ldots,\frac{W_n(H_N,\cdot)}{\sqrt{n}}
\right)
\end{equation}
converges as $n\to\infty$ in the sense of finite-dimensional distributions to
\begin{eqnarray*}
\left(
\sigma_q\,\sqrt{q!}\,B_q,\ldots,\sigma_N\,\sqrt{N!}\,B_N
\right),
\end{eqnarray*}
where $\sigma_s^2:=\sum_{k\in\Z}\rho(k)^s$ ($s=q,\ldots,N$), and $B_q,\ldots,B_N$ are
independent standard Brownian motions.
(The fact that $\sum_{k\in\Z}\rho(k)^s\geq 0$ is part of the conclusion.)
On the other hand, using (\ref{lien-tch}) as well as its Gaussian counterpart
(where one replaces, in (\ref{lien-tch}), the Tchebycheff polynomial $U_q$ and
the free Brownian motion $S$ by the Hermite polynomial $H_q$ and the standard Brownian
motion $B$ respectively; see, e.g., \cite[Theorem 2.7.7]{ivangiobook}),
we get, for any $s=q,\ldots,N$, that
\[
V_n(U_s,t)=I^{S}_s\left(\sum_{k=1}^{[nt]}e_k^{\otimes s}\right)
\quad\mbox{and}\quad
W_n(H_s,t)=I^B_s\left(\sum_{k=1}^{[nt]}e_k^{\otimes s}\right)
,
\]
where the sequence $\{e_k:\,k\in\Z\}$ is as in Section \ref{sec:semfree} and
$I^B_s(\cdot)$ stands for the multiple Wiener integral of order $s$
with respect to $B$. We observe that the kernel $\sum_{k=1}^{[nt]}
e_k^{\otimes s}$ is a symmetric
function of $L^2(\R^s)$.
Therefore, according to Proposition \ref{transfer} (transfer principle),
we deduce that the non-commutative counterpart of (\ref{h_nu}) holds as well, that is, we have that
\begin{eqnarray*}
\left(
\frac{V_n(U_q,\cdot)}{\sqrt{n}},\ldots,\frac{V_n(U_N,\cdot)}{\sqrt{n}}
\right),
\end{eqnarray*}
converges as $n\to\infty$ in the sense of finite-dimensional distributions to
\begin{eqnarray*}
\left(
\sigma_q\,S_q,\ldots,\sigma_N\,S_N
\right),
\end{eqnarray*}
where $S_q,\ldots,S_N$ denote freely independent free Brownian motions.
The desired conclusion (\ref{tcl1}) follows then as a consequence of this latter convergence, together with the decomposition
(\ref{decompo})
of $Q$ and
 the identity in law (see Section \ref{semi-sec}):
\[
a_q \sigma_qS_q + \ldots + a_N \sigma_NS_N \,\overset{\rm law}{=}\,\sqrt{a_q^2\sigma_q^2+\ldots+a_N^2\sigma_N^2}\times S.
\]
\qed

\subsection{Proof of the non-central limit theorem (\ref{tcl2})}

The proof is more delicate than the one for (\ref{tcl1}). This is
because the limit in the usual probability setting is not Gaussian, since it
is given by a multiple Wiener integral of order greater than 1.
Therefore, we cannot use the transfer principle as in (\ref{tcl1}). We need to focus
on the detailed structure of
\begin{equation}\label{e:VUU}
V_n(Q,t)=a_q\sum_{k=1}^{[nt]}U_q(X_k)+\sum_{s=q+1}^N a_s \sum_{k=1}^{[nt]}U_s(X_k).
\end{equation}
where we have again chosen $N$ large enough so that $a_s=0$ for all $s\geq N$.
The idea of the proof is to  show that, after normalization, the second term in (\ref{e:VUU}) is
asymptotically negligible, so that it is sufficient to focus on the first term
\[
G_n(t)=\sum_{k=1}^{[nt]}U_q(X_k).
\]
We can therefore apply Proposition \ref{l:key} which provides, in (\ref{kthmomentdiscret}),
an expression for $\ff(G_n(t_1)\ldots G_n(t_p))$ involving multiple sums. We show that the diagonals
in these multiple sums can be excluded (this is step 1 below). We express the remainder as integrals (step 2 below)
and apply the dominated convergence theorem to obtain the expression (\ref{kthmoment2})
which caracterizes the $q$th non-commutative Tchebycheff process $R_{H,q}$.

More specifically, fix
$\e\in(0,1]$. By virtue of (\ref{rho}), there exists an integer $M>0$ large enough so that, for all $j>M$,
\[
|\rho(j)|=j^{-D}L(j)\leq \e\leq 1.
\]
For any real $t>0$ and any integer $s$ larger than or equal to $q+1$, we can write
\begin{eqnarray}
&&\varphi\left[
\left(\frac{1}{n^{1-qD/2}L(n)^{q/2}}\sum_{k=1}^{[nt]} U_s(X_k)\right)^2
\right]
=\frac{1}{n^{2-qD}L(n)^q}\sum_{k,l=1}^{[nt]}\rho(k-l)^s\notag\\
&\leq&
\frac{1}{n^{2-qD}L(n)^q}\left(
\sum_{\substack{k,l=1\\|k-l|\leq M}}^{[nt]} 1
+2\e \sum_{\substack{k,l=1\\k> l+M}}^{[nt]} |\rho(k-l)|^q
\right)\notag\\
&\leq&\frac{t}{n^{1-qD}L(n)^q}\left((2M+1)+2\e\sum_{j=1}^{[nt]}j^{-qD}L(j)^q \right).\label{kata}
\end{eqnarray}
Since $qD<1$, we have that $n^{qD-1}L(n)^{-q}\to 0$ (to see this, use (\ref{potter}) with $1/L$ instead of $L$,
$1/L$ being slowly varying as well).
The following lemma is useful at this stage.
\begin{lemma}
When $t>0$ is fixed, we have
\begin{equation}\label{kara2}
\sum_{j=1}^{[nt]}j^{-qD}L(j)^q \sim \frac{[nt]^{1-qD} L([nt])^q}{1-qD}
\quad\mbox{as $n\to\infty$}.
\end{equation}
\end{lemma}
\noindent
{\it Proof}.
Although this is a somehow standard result in the theory of regular variation (Karamata's type theorem),
we prove (\ref{kara2}) for sake of completeness.
First, observe that
\[
\frac{1}{L([nt])^{q}[nt]^{1-qD}}\sum_{j=1}^{[nt]}j^{-qD}L(j)^q = \int_0^1 l_n(x)dx,
\]
where
\[
l_n(x)=\sum_{j=1}^{[nt]}\left(\frac{L(j)}{L([nt])}\right)^q\,\left(\frac{j}{[nt]}\right)^{-qD}
{\bf 1}_{[\frac{j-1}{[nt]}\frac{j}{[nt]})}(x).
\]
Since $L(j)/L([nt])=L([nt]\times(j/[nt]))/L([nt])\to 1$ for fixed $j/[nt]$ as $n\to\infty$,
one has $l_n(x)\to l_\infty(x)$ for $x\in (0,1)$, where $l_\infty(x)=x^{-qD}$.
By choosing a small enough $\delta>0$ so that $q(D+\delta)<1$ (this is possible because $qD<1$)
and $L(j)/L([nt])\leq C(j/[nt])^{-\delta}$ (this is possible thanks to (\ref{potter})), we get that
\[
|l_n(x)|\leq C\sum_{j=1}^{[nt]}\left(\frac{j}{[nt]}\right)^{-q(D+\delta)}
{\bf 1}_{[\frac{j-1}{[nt]}\frac{j}{[nt]})}(x)\leq Cx^{-(D+\delta)q}\quad\mbox{for all $x\in(0,1)$.}
\]
The function in the bound is integrable on $(0,1)$. Hence, the dominated convergence
theorem yields
\[
\frac{1}{L([nt])^{q}[nt]^{1-qD}}\sum_{j=1}^{[nt]}j^{-qD}L(j)^q\to \int_0^1 l_\infty(x)dx =\int_0^1 x^{-qD}dx=\frac{1}{1-qD},
\]
which is equivalent to (\ref{kara2}).
\qed

\bigskip

Let us go back to the proof of the non-central limit theorem. We have (\ref{kata}). But, since $L([nt])/L(n)\to 1$ ($t$ is fixed), we actually get that
\begin{equation}\label{kara}
\sum_{j=1}^{[nt]}j^{-qD}L(j)^q
\sim \frac{n^{1-qD}t^{1-qd} L(n)^q}{1-qD}\quad\mbox{as $n\to\infty$},
\end{equation}
so that, by combining (\ref{kara}) and (\ref{kata}),
\[
\limsup_{n\to\infty}\varphi\left[
\left(\frac{1}{n^{1-qD/2}L(n)^{q/2}}\sum_{k=1}^{[nt]} U_s(X_k)\right)^2
\right]\leq\frac{t^{2-qD}\,\e}{1-qD}.
\]
Since $\e>0$ is arbitrary, this implies that
\[
\frac{1}{n^{1-qD/2}L(n)^{q/2}}\sum_{k=1}^{[nt]} U_s(X_k)\overset{\rm L^2}{\to} 0\quad\mbox{as $n\to\infty$}
\]
for all $s\geq q+1$.
As a consequence, in the rest of the proof we can assume without loss of generality that $Q=a_qU_q$.
Thus,
set \[
F_n(t)=\frac{a_q}{n^{1-qD/2}L(n)^{q/2}}\sum_{k=1}^{[nt]}U_q(X_k) = \frac{a_q}{n^{1-qD/2}L(n)^{q/2}}\, G_n(t),
\quad t\geq 0,
\]
where $G_n$ is given by (\ref{gnn}).
Using (\ref{kthmomentdiscret}), we have that
\begin{equation}\label{split}
\varphi(F_n(t_1)\ldots F_n(t_p))=
\frac{(a_q)^p}{n^{p-pqD/2}L(n)^{pq/2}}
\sum_{k_1=1}^{[nt_1]}\ldots\sum_{k_p=1}^{[nt_p]}\,\,\,
\sum_{(r_1,\ldots,r_{p-1})\in B({\bf q}_p)}\,\,
\prod_{1\leq i<j\leq p}\rho(k_i-k_j)^{\alpha_{ij}({\bf r})},
\end{equation}
where ${\bf q}_p=(q_1,\ldots,q_p)=(q,\ldots,q)$.

To obtain the limit of (\ref{split}) as $n\to\infty$ and thus to conclude the proof of (\ref{tcl2}), we proceed in five steps.\\

{\it Step 1 (Determination of the main term)}. We split the sum $\sum_{k_1=1}^{[nt_1]}\ldots\sum_{k_p=1}^{[nt_p]}$ in the right-hand side of (\ref{split})
into
\begin{equation}\label{bos}
\sum_{\substack{k_1=1,\ldots,[nt_1]\\ ...\\k_p=1,\ldots,[nt_p]\\ \forall i\neq j:\, |k_i-k_j|\geq 3}}
\quad +\quad
\sum_{\substack{k_1=1,\ldots,[nt_1]\\ ...\\k_p=1,\ldots,[nt_p]\\ \exists i\neq j:\, |k_i-k_j|\leq 2}},
\end{equation}
and we show that the second sum in (\ref{bos}) is asymptotically negligible as $n\to\infty$.

Up to reordering, it is enough to show that, for any ${\bf r}=(r_1,\ldots,r_{p-1})\in B_p$,
\begin{equation}\label{boston}
R_n:=\frac{1}{n^{p-pqD/2}L(n)^{pq/2}}
\sum_{\substack{k_1=1,\ldots,[nt_1]\\ ...\\k_p=1,\ldots,[nt_p]\\ |k_{p-1}-k_p|\leq 2}}
\,\,\,
\prod_{1\leq i<j\leq p}\rho(k_i-k_j)^{\alpha_{ij}({\bf r})}
\end{equation}
tends to zero as $n\to\infty$. In (\ref{boston}), let us bound $|\rho(k_i-k_p)|$ by $1$
when $i\in\{1,\ldots,p-1\}$.
We get that
\begin{eqnarray}
\notag
R_n\leq \frac{5}{n^{p-pqD/2}L(n)^{pq/2}}
\sum_{k_1=1}^{[nt_1]}\ldots\sum_{k_{p-1}=1}^{[nt_{p-1}]}
\,\,\,
\prod_{1\leq i<j\leq p-1}\rho(k_i-k_j)^{\alpha_{ij}({\bf r})}.\\
\label{boston2}
\end{eqnarray}
Going back to Definition \ref{def:alpha}, there is $\widehat{{\bf q}}_{p-1}=(\widehat{q}_1,\ldots,\widehat{q}_{p-1})$ and
$\widehat{{\bf r}}=(\widehat{r}_1,\ldots,\widehat{r}_{p-2})\in B(\widehat{{\bf q}}_{p-1})$
such that $\alpha_{ij}(\widehat{{\bf r}})=\alpha_{ij}({\bf r})$ for all $i,j=1,\ldots,p-1$.
This is because the connexions between the remaining $p-1$ blocks are unchanged.
Moreover, since we remove the $q$ connexions associated to the block $p$  (this involves $2q$ dots, see Figure 3), we have
\begin{equation}\label{sommedesqchapeaux}
\widehat{q}_1+\ldots+\widehat{q}_{p-1}=(p-2)q.
\end{equation}
Hence, using Lemma \ref{quivaienb}
with $T=\Z$,
$\nu_i=\sum_{k=1}^{[nt_i]}\delta_k$ ($\delta_k$ being the Dirac mass at $k$)
and $e(k,x)=e_k(x)$ as in (\ref{identityquivabien}),
we get that
\begin{eqnarray}
\sum_{k_1=1}^{[nt_1]}\ldots\sum_{k_{p-1}=1}^{[nt_{p-1}]}
\,\,\,
\prod_{1\leq i<j\leq p-1}\rho(k_i-k_j)^{\alpha_{ij}({\bf r})}=
(\ldots((\widehat{g}_{1,n}(t_1,\cdot)\overset{\widehat{r}_1}\frown \widehat{g}_{2,n}(t_2,\cdot))\overset{\widehat{r}_2}
\frown\ldots)
\overset{\widehat{r}_{p-2}}\frown \widehat{g}_{p-1,n}(t_{p-1},\cdot)\notag\\
\label{nyc}
\end{eqnarray}
where, for any $i=1,\ldots,p-1$,
\[
\widehat{g}_{i,n}(t_i,x_1,\ldots,x_{\widehat{q}_i})=\sum_{k=1}^{[nt_i]}e_k(x_1)\ldots e_k(x_{\widehat{q}_i}).
\]
Iterative applications of (\ref{cauchy}) in (\ref{nyc}) lead to
\[
\sum_{k_1=1}^{[nt_1]}\ldots\sum_{k_{p-1}=1}^{[nt_{p-1}]}
\,\,\,
\prod_{1\leq i<j\leq p-1}|\rho(k_i-k_j)|^{\alpha_{ij}({\bf r})}
\leq
\|\widehat{g}_{1,n}(t_1,\cdot)\|_{L^2(\R^{\widehat{q}_1})}
\ldots
\|\widehat{g}_{p-1,n}(t_{p-1},\cdot)\|_{L^2(\R^{\widehat{q}_{p-2}})}.
\]
But, for any $i=1,\ldots,p-1$,
\begin{eqnarray*}
\|\widehat{g}_{i,n}(t_i,\cdot)\|_{L^2(\R^{\widehat{q}_i})}^2
&=&\sum_{k,l=1}^{[nt_i]}\rho(k-l)^{\widehat{q}_i}
=\sum_{|j|<[nt_i]}\rho(j)^{\widehat{q}_i}\big([nt_i]-|j|\big)\\
&=& [nt_i]+2\sum_{1\leq j<[nt_i]} j^{-\widehat{q}_i D}L(j)^{\widehat{q}_i}\big([nt_i]-j\big)\\
&\leq& [nt_i]\left(1+
2\sum_{1\leq j<[nt_i]} j^{-\widehat{q}_i D}L(j)^{\widehat{q}_i}
\right)
\leq Cn^{2-\widehat{q}_iD}L(n)^{\widehat{q}_i},
\end{eqnarray*}
for some $C>0$,
and where the last inequality holds because of (\ref{kara}).
We deduce
\begin{eqnarray*}
\sum_{k_1=1}^{[nt_1]}\ldots\sum_{k_{p-1}=1}^{[nt_{p-1}]}
\,\,\,
\prod_{1\leq i<j\leq p-1}|\rho(k_i-k_j)|^{\alpha_{ij}({\bf r})}
&\leq&
Cn^{p-1-(\widehat{q}_1+\ldots+\widehat{q}_{p-1})D/2}L(n)^{(\widehat{q}_1+\ldots+\widehat{q}_{p-1})/2}\\
&=&
Cn^{p-1-(p-2)qD/2}L(n)^{(p-2)q/2},
\end{eqnarray*}
by (\ref{sommedesqchapeaux}).
By putting all these bounds together, and because $qD<1$,we conclude from (\ref{boston2}) that $R_n$ given by (\ref{boston}) tends to zero as
$n\to\infty$.
\medskip

{\it Step 2 (Expressing sums as integrals)}.
We now consider the first term in (\ref{bos}) and express it as an integral, so to apply
the dominated convergence theorem.

From Definition \ref{def:alpha}, we deduce immediately that $\sum_{1\leq i<j\leq p}\alpha_{ij}({\bf r})=pq/2$
when ${\bf r}=(r_1,\ldots,r_{p-1})\in B({\bf q}_p)$.
This fact, combined with the specific form (\ref{rho}) of $\rho$, yields
\begin{eqnarray*}
&&\frac{(a_q)^p}{n^{p-pqD/2}L(n)^{pq/2}}
\sum_{\substack{k_1=1,\ldots,[nt_1]\\ ...\\k_p=1,\ldots,[nt_p]\\ \forall i\neq j:\, |k_i-k_j|\geq 3}}
\,\,\,
\sum_{(r_1,\ldots,r_{p-1})\in B({\bf q}_p)}\,\,
\prod_{1\leq i<j\leq p}\rho(k_i-k_j)^{\alpha_{ij}({\bf r})}\\
&=&
\frac{(a_q)^p}{n^{p}}
\sum_{\substack{k_1=1,\ldots,[nt_1]\\ ...\\k_p=1,\ldots,[nt_p]\\ \forall i\neq j:\, |k_i-k_j|\geq 3}}
\,\,\,\sum_{(r_1,\ldots,r_{p-1})\in B({\bf q}_p)}\,\,
\prod_{1\leq i<j\leq p}
\left[\left|\frac{k_i-k_j}{n}
\right|^{-D\alpha_{ij}({\bf r})}\left(\frac{L(|k_i-k_j|)}{L(n)}\right)^{\alpha_{ij}({\bf r})}\right]\\
&=&(a_q)^p\,\,\,\sum_{(r_1,\ldots,r_{p-1})\in B({\bf q}_p)}\,\,\int_{0}^{\infty}\ldots \int_{0}^{\infty}l_{n,{\bf r}}(s_1,\ldots,s_p)ds_1\ldots ds_p,
\end{eqnarray*}
where
\begin{eqnarray*}
l_{n,{\bf r}}(s_1,\ldots,s_p)&=&
\sum_{\substack{k_1=1,\ldots,[nt_1]\\ ...\\k_p=1,\ldots,[nt_p]\\ \forall i\neq j:\, |k_i-k_j|\geq 3}}
\prod_{1\leq i<j\leq p}\left[\left|\frac{k_i-k_j}{n}
\right|^{-D\alpha_{ij}({\bf r})}\left(\frac{L(|k_i-k_j|)}{L(n)}\right)^{\alpha_{ij}({\bf r})}\right]\\
&&\hskip6cm\times
{\bf 1}_{[\frac{k_1-1}{n},\frac{k_1}n)}(s_1)\ldots
{\bf 1}_{[\frac{k_p-1}{n},\frac{k_p}n)}(s_p).
\end{eqnarray*}

\medskip

{\it Step 3 (Pointwise convergence)}.
We show the pointwise convergence of $l_{n,{\bf r}}$.
Since, for fixed $|k_i-k_j|/n$ and as $n\to\infty$, one
has
 \[
L(|k_i-k_j|)/L(n)=L(n\times|k_i-k_j|/n)/L(n)\to 1,
\]
one deduces that
 $l_{n,{\bf r}}(s_1,\ldots,s_p)\to l_{\infty,{\bf r}}(s_1,\ldots,s_p)$ for any $s_1,\ldots,s_p\in \R_+$, where
\[
l_{\infty,{\bf r}}(s_1,\ldots,s_p)=
{\bf 1}_{[0,t_1]}(s_1)\ldots {\bf 1}_{[0,t_p]}(s_p)
\prod_{1\leq i<j\leq p}\left|s_i-s_j
\right|^{-D\alpha_{ij}({\bf r})}.
\]

\medskip

{\it Step 4 (Domination)}.
We show that $l_{n,{\bf r}}$ is dominated by an integrable function.
If $k_i-k_j\geq 3$ (the case where $k_j-k_i\geq 3$ is similar by symmetry),
$s_i\in [\frac{k_i-1}{n},\frac{k_i}n)$ and
$s_j\in [\frac{k_j-1}{n},\frac{k_j}n)$, then
\[
\frac{3}{n}\leq \frac{k_i-k_j}{n}\leq s_i+\frac{1}{n}-s_j,
\]
so that $s_i-s_j\geq \frac{2}{n}$, implying in turn
\begin{equation}\label{kikj}
\frac{k_i-k_j}{n}\geq s_i-s_j-\frac{1}{n}\geq \frac{s_i-s_j}{2}.
\end{equation}
Since $qD<1$, choose a small enough $\delta$ so that $q(D+\delta)<1$ and
$L(k_i-k_j)/L(n)\leq C((k_i-k_j)/n)^{-\delta}$
(this is possible thanks to (\ref{potter})).
We get that
\begin{eqnarray*}
|l_{n,{\bf r}}(s_1,\ldots,s_p)|&\leq& C
\sum_{\substack{k_1=1,\ldots,[nt_1]\\ ...\\k_p=1,\ldots,[nt_p]\\ \forall i\neq j:\, |k_i-k_j|\geq 3}}
\prod_{1\leq i<j\leq p}\left|\frac{k_i-k_j}{n}
\right|^{-(D+\delta)\alpha_{ij}({\bf r})}
{\bf 1}_{[\frac{k_1-1}{n},\frac{k_1}n)}(s_1)\ldots
{\bf 1}_{[\frac{k_p-1}{n},\frac{k_p}n)}(s_p)\\
&\leq&C\,2^{(D+\delta)pq/2}{\bf 1}_{[0,t_1]}(s_1)\ldots {\bf 1}_{[0,t_p]}(s_p)
\prod_{1\leq i<j\leq p}\left|s_i-s_j
\right|^{-(D+\delta)\alpha_{ij}({\bf r})},
\end{eqnarray*}
by (\ref{kikj}).
The function in the bound is integrable on $\R_+^p$. Indeed, for any fixed
${\bf r}=(r_1,\ldots,r_{p-1})\in B({\bf q}_p)$, we can write, thanks to (\ref{kthmomentfinal2}) and with
$f_{1-q(D+\delta)/2,q}$ given by (\ref{e:kernel}),
\begin{eqnarray*}
&&\int_0^{t_1}ds_1\ldots \int_0^{t_p}ds_p\prod_{1\leq i<j\leq p}\left|s_i-s_j
\right|^{-(D+\delta)\alpha_{ij}({\bf r})}\\
&=&C\,
(\ldots((f_{1-q(D+\delta)/2,q}(t_1,\cdot)\overset{r_1}\frown f_{1-q(D+\delta)/2,q}(t_2,\cdot))\overset{r_2}\frown
\ldots)
\overset{r_{p-1}}\frown f_{1-q(D+\delta)/2,q}(t_p,\cdot),
\end{eqnarray*}
where $C>0$,
so that by an iterative use of (\ref{cauchy}),
\begin{eqnarray*}
&&\int_0^{t_1}ds_1\ldots \int_0^{t_p}ds_p\prod_{1\leq i<j\leq p}\left|s_i-s_j
\right|^{-(D+\delta)\alpha_{ij}({\bf r})}\\
&\leq&
\|f_{1-q(D+\delta)/2,q}(t_1,\cdot)\|_{L^2(\R^q)}\ldots \|f_{1-q(D+\delta)/2,q}(t_p,\cdot)\|_{L^2(\R^q)}
=(t_1\ldots t_p)^{1-q(D+\delta)/2}<\infty,
\end{eqnarray*}
where the last equality holds because of (\ref{335}).\\

{\it Step 5 (Dominated convergence)}. By combining the results of Steps 2 to 4, we obtain that the dominated convergence
theorem applies and yields
\begin{eqnarray*}
\varphi(F_n(t_1)\ldots F_n(t_p))\to
(a_q)^p\int_0^{t_1} ds_1\ldots \int_0^{t_p} ds_p\,\,\,
\sum_{(r_1,\ldots,r_{p-1})\in B({\bf q}_p)}\,\,
\prod_{1\leq i<j\leq p}|s_i-s_j|^{-D\alpha_{ij}({\bf r})}.
\end{eqnarray*}
We recognize that, up to a multiplicative constant, this is the quantity in (\ref{kthmoment2}) with $H=1-qD/2$.
More precisely, we have
\[
\varphi(F_n(t_1)\ldots F_n(t_p))\to
\varphi\left(\frac{a_q}{\sqrt{H(2H-1)}}\, R_{H,q}(t_1)\times\ldots\times \frac{a_q}{\sqrt{H(2H-1)}}\, R_{H,q}(t_p)\right),
\]
which concludes the proof of (\ref{tcl2}).
\qed

\bigskip

\noindent
{\bf Acknowledgments}.
We would like to thank two anonymous referees for their careful
reading of the manuscript and for their valuable suggestions and remarks. Also,
I. Nourdin would like to warmly thank M. S. Taqqu
for his hospitality during his stay at Boston University in October 2011, where part of this research was carried out.

\end{document}